\documentclass[12pt]{article}
\usepackage[margin=1in]{geometry}
\DeclareMathAlphabet\mathbfcal{OMS}{cmsy}{b}{n}
\usepackage[utf8]{inputenc}
\usepackage[dvipsnames]{xcolor}
\usepackage{enumerate}
\usepackage[utf8]{inputenc}
\usepackage{amsfonts}
\usepackage{amsthm}
\usepackage{amsmath}
\usepackage{array}
\usepackage{lscape}
\usepackage{mathrsfs}
\usepackage{graphicx}
\usepackage{tikz-cd}
\usepackage{calc}
\usepackage{float}
\usepackage{stmaryrd}
\usepackage{multirow}
\usepackage{lipsum}
\usepackage{amssymb}
\usepackage{mathrsfs}
\usepackage{mathtools}
\usepackage{mathdots}
\usepackage{fancyhdr}
\usepackage{tikz}
\usepackage{faktor}
\usepackage{setspace}
\usetikzlibrary{decorations.markings}
\newtheorem{theorem}{Theorem}
\theoremstyle{plain}
\usepackage{sectsty}
\usepackage{kpfonts}
\newtheorem{corollary}{Corollary}

\newtheorem{example}{Example}
\newtheorem{lemma}{Lemma}

\newtheorem{proposition}{Proposition}
\newtheorem{remark}{Remark}

\newtheorem{question}{Question}

\newtheorem*{proposition*}{Proposition}

\newcommand{\GL}{{\mbox{GL}}}
\newcommand{\Hom}{{\mbox{Hom}}}

\newcommand{\End}{{\mbox{End}}}

\newcommand{\vol}{{\mbox{vol}}}
\newcommand{\Fr}{{\mbox{Fr}}}

\newtheorem*{observation*}{Observation}
\newtheorem*{theorem*}{Theorem}
\newtheorem*{claim*}{Claim}  
\usepackage{blindtext}
\usepackage{fancyhdr}
\title{A curvature obstruction to integrability}
\author{Gabriella Clemente}
\date{}
\allowdisplaybreaks
\makeatletter
\renewcommand\tableofcontents{%
    \@starttoc{toc}%
}
\begin{document}
\maketitle

\begin{abstract}
The classical theory of $G$-structures, which include almost-complex structures, explains the relationship between the curvature of compatible connections and integrability. This note is an effort to understand how the curvature of Riemannian metrics can obstruct the integrability of almost-complex structures. It is shown that certain special complex structures cannot coexist with non-flat constant curvature metrics, and a formal variational realization of these structures is provided. The approach followed here is direct, meaning that it bypasses the classical theory. The idea is to find obstruction equations for the integrability of almost-complex structures by way of Nijenhuis tensor derivatives. These new equations involve the curvature of a torsion-free connection, and reveal the interplay between the almost-complex and Riemannian geometries. Curvature scalars to detect non-complexity in the compact case then arise in a natural way.    
\end{abstract}

\section{Introduction}
The motivating question of this note is 

\begin{question}\label{IQ}
Do (compact) almost-complex but not complex manifolds of real dimension at least $6$ exist?
\end{question}

Here it will be shown that in high dimensions, non-flat, constant curvature Riemannian metrics obstruct the existence of certain special complex structures (Theorem \ref{TV1}).

Question \ref{IQ} is a differentio-geometrical obstruction theory problem. The classical obstruction theory for $G$-structures, of which almost-complex structures are an example, explains the link between integrability and the curvature of certain linear connections. This will be outlined below, and most of the details can be found in \cite{Cap} and \cite{Cat}. However, the approach used in the present note was developped completely independently from this theory. Let $G \leq \GL_n(\mathbb{R})$ be a Lie group. Recall that a $G$-structure on a manifold $M^n$ is a $G$-invariant submanifold $P \subseteq \Fr(M) \xrightarrow{\pi} M$ of the frame bundle of $M$ such that $\pi|_P:P \to M$ is a principal $G$-bundle. And $P$ is integrable if it is locally isomorphic, in the $G$-structural sense, to the flat $G$-structure $\mathbb{R}^n \times G$ on $\mathbb{R}^n.$ In a nutshell, the intrinsic torsions $T^k_{intr} (P),$ for $k \geq 1$ an integer, are Spencer cohomology classes obstructing the integrability of $P.$ Let $\nabla$ be a linear connection on $M$ that is compatible with $P.$ Namely, the associated connection $1$-form $\omega^{\nabla} \in \Omega^1 \big(\Fr(M),\mathfrak{gl}_n(\mathbb{R})\big)$ satisfies $\omega^{\nabla}|_P \in \Omega^1 \big(P,\mathfrak{g}\big),$ where $\mathfrak{g}$ is the Lie algebra of $G.$ The first $2$ intrinsic torsions have the following geometrical interpretation: $T^1_{intr} (P)=[T^{\nabla}]$ and $T^2_{intr} (P)=[R^{\nabla}]$; i.e.\ the Spencer classes $T^i_{intr}(P),$ $1 \leq i \leq 2,$ are represented by the torsion, respectively the curvature of the connection $\nabla.$ For example, an almost-complex structure on $M^{2n}$ is a $\GL_n(\mathbb{C})$-structure, and the torsion of any compatible connection coincides with the Nijenhuis tensor (see section \ref{int} for its defining formula). Analogously, a $G$-structure $P$ equipped with the connection $1$-form of any compatible connection $\nabla,$ $(\pi:P \to M^n, \omega^{\nabla}),$ can be viewed as a Cartan geometry of type $(\mathbb{R}^n \rtimes G, G).$ The Cartan connection is $\Omega_{\omega^{\nabla}}:=\theta+\omega^{\nabla} \in \Omega^1\big(P,\mathbb{R}^n \oplus \mathfrak{g}\big),$ where $\theta \in \Omega^1 \big(P, \mathbb{R}^n \big)$ is the soldering form of $P.$ The curvature form $d\Omega_{\omega^{\nabla}}+\frac{1}{2}[\Omega_{\omega^{\nabla}}, \Omega_{\omega^{\nabla}}]$ of the Cartan connection contains both the torsion and curvature forms of $\omega^{\nabla}.$ The latter two correspond to the torsion and curvature of the linear connection $\nabla.$ 

The almost-complex and Riemannian geometries of a manifold exists independently from each another, and are essentially different kinds of geometries. But they can be seen to interact. For instance, twistor spaces have been used as an extrinsic tool for understanding how local conformal flatness relates to orthogonal complex structures (see \cite{Salamon}, and the sources cited therein). The goal of this note is to probe the almost-complex geometry with Riemannian metrics of prescribed curvature. The metrics in focus are the non-zero constant sectional curvature ones. It should be possible to unify this point of view with the $G$-structure perspective. But that will be subject of future investigations.  

Theorem \ref{TV1} gives rise to a more tractable version of Question \ref{IQ}: can a compact, almost-complex Riemannian manifold of positive constant sectional curvature be complex? A study of this question can result in insights about the (non-)integrability of almost-complex structures on the six-sphere $S^6.$ A negative answer would seem to imply that a high dimensional, compact, complex manifold admitting a Riemannian metric of everywhere positive sectional curvature must be ridged, in the sense that the curvature of any such metric must have some undulation, even if ever so slight. If $S^6$ were non-complex, one could then try to imagine its ridged dual, a ridged sphere, a compact, complex correlate manifold that is curvature-wise the closest one to the standard sphere.

The methodology of this note is to set up a calculus of vector bundle valued differential forms that allows for an intrinsic study of the (almost-)complex geometry of manifolds in general. The second section covers the algebraic structures on spaces of bundle forms that is needed for this calculus. The third section explains how Nijenhuis derivatives interface the metric and (almost-)complex geometries. The curvature obstruction equations (see Propoition \ref{PPP1} and Remark \ref{REM1}) are well-suited for ruling out complex structures that satisfy additional algebraic relations, such as the previously mentioned special complex structures. They seem to provide only point-wise constraints. The global picture must take into account the topology of the manifold, and can be captured with numerical obstructions obtained from Nijenhuis tensor derivatives (curvature scalars). The last 2 sections show how special hypothetical complex structures interact with sectional curvature, how to variationally realize them, and how non-existence of a complex structure could, in principle, be detected via curvature scalars (Proposition \ref{invC}).

\section{Algebraic tools}\label{algsect}
The calculus of vector bundle valued forms that appears starting from the next section stands on an algebraic foundation that is expanded on here. 

Let $M^n$ be a smooth manifold, and $\Omega^{\bullet} (M, T_M)=\bigoplus_{k=0}^n \Omega^k (M, T_M)$ be the space of tangent bundle valued differential forms. Equip the spaces \[\Omega^{\bullet} \big(M, \End_{\mathbb{R}}(T_M)\big)=\bigoplus_{k=0}^n \Omega^k \big(M, \End_{\mathbb{R}}(T_M)\big), \mbox{ and } \Omega^{\bullet} \big(M, \bigwedge^{\bullet} {T_M}\big)=\bigoplus_{k=0}^{2n} \bigoplus_{p+q=k}  \Omega^p \big(M, \bigwedge^q {T_M}\big)\] of endomorphism, respectively polyvector valued forms with operations defined as follows. 

Consider the action \[\cdot:\End_{\mathbb{R}}(T_M) \times T^*_M \otimes T_M \to T^*_M \otimes T_M,\] \[(S, f\otimes e)\mapsto S\cdot (f\otimes e):=f\otimes S(e).\] The product in $\Omega^{\bullet} \big(M, \End_{\mathbb{R}}(T_M)\big)$ is given as follows: for any $\alpha \in \Omega^k \big(M, \End_{\mathbb{R}}(T_M)\big), \beta \in \Omega^l \big(M, \End_{\mathbb{R}}(T_M)\big),$ $\alpha \wedge \beta \in \Omega^{k+l} \big(M, \End_{\mathbb{R}}(T_M)\big),$ and 
\begin{equation*}
\begin{split}
(\alpha \wedge \beta)(X_1,\dots,X_{k+l})=\frac{1}{k!l!}\sum_{\sigma \in S_{k+l}} sign(\sigma) \alpha(X_{\sigma(1)},\dots,X_{\sigma(k)})\cdot \beta(X_{\sigma(k+1)},\dots,X_{\sigma(k+l)}),
\end{split}
\end{equation*}
where indeed \[\alpha(X_{\sigma(1)},\dots,X_{\sigma(k)}), \beta(X_{\sigma(k+1)},\dots,X_{\sigma(k+l)}) \in \End_{\mathbb{R}}(T_M).\]

The product in $\Omega^{\bullet} \big(M, \bigwedge^{\bullet} {T_M}\big)$ is defined with a normalization factor but without any twisting. Namely, for any $\gamma \in  \Omega^i \big(M, \bigwedge^j {T_M}\big), \theta \in \Omega^{k} \big(M, \bigwedge^{l} {T_M}\big),$ $\gamma \wedge \theta \in \Omega^{i+k} \big(M, \bigwedge^{j+l} {T_M}\big),$ and

\begin{equation*}
\begin{split}
(\gamma \wedge \theta)(X_1,\dots,X_{i+k})=\frac{1}{2}\frac{1}{i!k!}\sum_{\sigma \in S_{i+k}} sign(\sigma) \gamma(X_{\sigma(1)},\dots,X_{\sigma(i)})\wedge \theta(X_{\sigma(i+1)},\dots,X_{\sigma(i+k)}),
\end{split}
\end{equation*}
where \[\gamma(X_{\sigma(1)},\dots,X_{\sigma(i)})\in \bigwedge^j {T_M}, \mbox{ and }\theta(X_{\sigma(i+1)},\dots,X_{\sigma(i+k)})\in \bigwedge^{l} {T_M}.\] 

Thus, $\Omega^{\bullet} \big(M, \End_{\mathbb{R}}(T_M)\big),$ and $\Omega^{\bullet} \big(M, \bigwedge^{\bullet} {T_M}\big)$ can be regarded as graded algebras with these products.

Note that if $\gamma \in \Omega^i(M,T_M),$ and $\theta \in \Omega^k(M,T_M),$ then $\gamma \wedge \theta = (-1)^{ik+1} \theta \wedge \gamma$ so that $\gamma, \theta$ anti-commute iff one of $i$ and $k$ is even. 

The tangent bundle forms $\Omega^{\bullet} (M, T_M)$ are acted on the left by $\Omega^{\bullet} \big(M, \End_{\mathbb{R}}(T_M)\big),$ and on the right by $\Omega^{\bullet} \big(M, \bigwedge^{\bullet} {T_M}\big).$ The left action combines the evaluation map \[ \End_{\mathbb{R}}(T_M) \times T_M \to T_M, \quad (S,y) \mapsto S(y),\] with the wedge product of forms: for any $\rho \in \Omega^s (M, T_M),$ and $\alpha \in \Omega^k \big(M, \End_{\mathbb{R}}(T_M)\big),$ $\alpha \wedge \rho \in \Omega^{k+s}(M, T_M),$ and 

\begin{equation*}
\begin{split}
(\alpha \wedge \rho)(X_1,\dots,X_{k+s})=\frac{1}{k!s!}\sum_{\sigma \in S_{k+s}} sign(\sigma) \alpha(X_{\sigma(1)},\dots,X_{\sigma(k)})\big(\rho(X_{\sigma(k+1)},\dots,X_{\sigma(k+s)})\big),
\end{split}
\end{equation*}
where \[\alpha(X_{\sigma(1)},\dots,X_{\sigma(k)})\in \End_{\mathbb{R}}(T_M), \mbox{ and }\rho(X_{\sigma(k+1)},\dots,X_{\sigma(k+s)}) \in T_M.\]
The right action makes use of the isomorphism \[\Hom_{\mathbb{R}}\Big(\bigwedge^j {T_M},T_M\Big)\simeq \bigwedge^j {T^*_M} \otimes T_M,\] the evaluation map 
\[\Hom_{\mathbb{R}}\Big(\bigwedge^j {T_M},T_M\Big) \times \bigwedge^j {T_M} \to T_M,\]
\[(C:=c\otimes z,x)\mapsto C(x)=(c\otimes z)(x)=c(x)z,\] and the wedge product of $\Omega^{\bullet}(M).$

Indeed, if $\gamma \in  \Omega^i \big(M, \bigwedge^j {T_M}\big),$ and $s \geq j,$ set 

\begin{equation*}
\begin{split}
(\rho \wedge \gamma)(X_1,\dots,X_{s-j+i})=\frac{1}{(s-j)!i!}\sum_{\sigma \in S_{s-j+i}} sign(\sigma) \rho(X_{\sigma(1)},\dots,X_{\sigma(s-j)},\cdot,\dots,\cdot)\big(\gamma(X_{\sigma(s-j+1)},\dots,X_{\sigma(s-j+i)})\big),
\end{split}
\end{equation*}
where \[\rho(X_{\sigma(1)},\dots,X_{\sigma(s-j)},\cdot,\dots,\cdot) \in \Hom_{\mathbb{R}}\Big(\bigwedge^j {T_M},T_M\Big), \mbox{ and }\gamma(X_{\sigma(s-j+1)},\dots,X_{\sigma(s-j+i)}) \in \bigwedge^j {T_M}.\] If $s<j,$ declare $\rho \wedge \gamma=0.$

For simplicity, the same notation is being used to express all products at play here. However, parentheses and context do away with ambiguities.

\begin{example}\label{ia1}
If $s \geq j=2,$ $\gamma=\gamma_{i-a} \wedge \gamma_a$ for $0 \leq a \leq i,$ then 
\begin{equation*}
\begin{split}
(\rho \wedge \gamma)(X_1,\dots,X_{s-2+i})&=\frac{1}{(s-2)!i!} \sum_{\sigma \in S_{s-2+i}} sign(\sigma) \rho(X_{\sigma(1)},\dots,X_{\sigma(s-2)},\cdot,\dots,\cdot)\times \\
&\big((\gamma_{i-a} \wedge \gamma_a)(X_{\sigma(s-1)},\dots,X_{\sigma(s-2+i)})\big)\\
&=\frac{1}{2(s-2)! i! (i-a)! a!} \sum_{\sigma \in S_{s-2+i}} \sum_{\theta \in S_i} sign(\sigma) sign(\theta) \times \\
& \rho \big(X_{\sigma(1)},\dots,X_{\sigma(s-2)},\gamma_{i-a}(X_{\theta(\sigma(s-1))},\dots,X_{\theta(\sigma(s-2+i-a))}),\gamma_a(X_{\theta(\sigma(s-1+i-a))},\dots,\\
&X_{\theta(\sigma(s-2+i))})\big)
\end{split}
\end{equation*}
Particularly, for $\alpha \in \Omega^1 (M, T_M),$ $\beta \in \Omega^2 (M, T_M),$ 
\[\big(\beta \wedge (\alpha \wedge \alpha)\big)(X_1,X_2)=\beta \big(\alpha(X_1),\alpha(X_2)\big).\] The computation in full detail:
\begin{equation*}
\begin{split}
\big(\beta \wedge (\alpha \wedge \alpha)\big)(X_1, X_2)&=\frac{1}{2} \sum_{\sigma \in S_2} sign(\sigma) \beta\big((\alpha \wedge \alpha)(X_{\sigma(1)}, X_{\sigma(2)})\big),
\end{split}
\end{equation*}
where
\begin{equation*}
\begin{split}
(\alpha \wedge \alpha)(X_{\sigma(1)}, X_{\sigma(2)})&=\frac{1}{2} \sum_{\sigma' \in S_2} sign(\sigma') \alpha(X_{\sigma'(\sigma(1))}) \wedge\alpha(X_{\sigma'(\sigma(2))})\\
&=\frac{1}{2}\big(\alpha(X_{\sigma(1)}) \wedge\alpha(X_{\sigma(2)})-\alpha(X_{\sigma(2)})\wedge \alpha(X_{\sigma(1)}) \big)\\
&=\frac{1}{2} \big(2\alpha(X_{\sigma(1)}) \wedge\alpha(X_{\sigma(2)})\big)\\
&=\alpha(X_{\sigma(1)}) \wedge\alpha(X_{\sigma(2)}),
\end{split}
\end{equation*}
so
\begin{equation*}
\begin{split}
\big(\beta \wedge (\alpha \wedge \alpha)\big)(X_1, X_2)&=\frac{1}{2} \sum_{\sigma \in S_2} sign(\sigma) \beta \big(\alpha(X_{\sigma(1)}) \wedge\alpha(X_{\sigma(2)})\big)\\
&=\frac{1}{2} \sum_{\sigma \in S_2} sign(\sigma) \beta(\alpha(X_{\sigma(1)}), \alpha(X_{\sigma(2)}))\\
&=\frac{1}{2}[\beta(\alpha(X_1), \alpha(X_2))-\beta(\alpha(X_2), \alpha(X_1)) ]\\
&=\beta(\alpha(X_1), \alpha(X_2)).
\end{split}
\end{equation*}
\end{example}

\section{Curvature and integrability}\label{int}
Let \[AC(M):=\{A \in \Omega^1 (M, T_M) \mid A \circ A=-Id\}\] be the space of almost-complex structures on $M.$ Unless otherwise stated, $M$ will be assumed to be an almost-complex manifold ($AC(M)\neq \emptyset$) of real even dimension $n \geq 2.$ Recall that $M$ is complex if it carries an $A \in AC(M)$ such that the Nijenhuis tensor of $A,$ \[N_A (\zeta, \eta)=[A(\zeta), A(\eta)]-A([A(\zeta),\eta]+[\zeta,A(\eta)])-[\zeta,\eta],\] vanishes for all vector fields $\zeta, \eta \in \mathfrak{X}(M)$ \cite{New}. In this case, $A$ is called an integrable almost-complex structure or a complex structure. Let \[C(M)=\{A \in AC(M) \mid N_A=0\}\] be the space of complex structures on $M.$

From this point on, $\nabla$ is a symmetric connection on $T_M,$ and $d^{\nabla}$ is the covariant exterior derivative associated to it, which at degree $k,$ is the map $d^{\nabla}:\Omega^k (M, T_M) \to \Omega^{k+1} (M, T_M),$ 
\begin{equation*}
\begin{split}
(d^{\nabla} \alpha)(\zeta_0,\dots,\zeta_k)&=\sum_{i=0}^k (-1)^i \nabla_{\zeta_i} \alpha(\zeta_0,\dots,\widehat{\zeta_i},\dots,\zeta_k)+\\
&\sum_{0 \leq i <j \leq k} (-1)^{i+j} \alpha ([\zeta_i,\zeta_j],\dots,\widehat{\zeta_i},\dots,\widehat{\zeta_j},\dots,\zeta_k)\\
&=\sum_{i=0}^k (-1)^i (\nabla_{\zeta_i} \alpha) (\zeta_0,\dots,\hat{\zeta_i},\dots,\zeta_k).
\end{split}
\end{equation*}
Also, $R^{\nabla}$ will denote here the curvature of the connection $\nabla$; for any $X,Y,Z,W \in \mathfrak{X}(M),$ \[R^{\nabla}(X,Y)Z=\nabla_X \nabla_Y Z-\nabla_Y \nabla_X Z-\nabla_{[X,Y]} Z. \] 
The bundle form language of the previous section has the capacity to articullate the integrability condition for almost-complex structures. 
\begin{lemma}\label{Q1}
Let $A \in AC(M).$ Then, $A$ is integrable iff \[d^{\nabla} A \wedge (A \wedge A)-d^{\nabla} A=0.\]
\end{lemma}

\begin{proof}
The Nijenhuis tensor of $A$ can be rewritten in terms of $d^{\nabla},$ using the condition $A \circ A=-Id$:

\begin{equation}\label{fog}
\begin{split}
A \circ \big(d^{\nabla} A(A(\zeta), A(\eta))-d^{\nabla} A(\zeta, \eta)\big)&=A \circ \Big[ \big((\nabla_{A(\zeta)} A)(A(\eta))-(\nabla_{A(\eta)} A)(A(\zeta))\big)\\
&-\big( (\nabla_{\zeta} A)(\eta)-(\nabla_{\eta} A)(\zeta)\big)\Big]\\
&=A \circ \Big[-A \big(\nabla_{A(\zeta)} A(\eta) -  \nabla_{A(\eta)} A(\zeta)\big)-\Big(\big(\nabla_{A(\zeta)} \eta- \nabla_{\eta} A(\zeta) \big)+\\
&\big(\nabla_{\zeta} A(\eta)- \nabla_{A(\eta)} \zeta \big)\Big)+A(\nabla_{\zeta} \eta-\nabla_{\eta} \zeta)\Big]\\
&=[A(\zeta), A(\eta)]-A\big([A(\zeta), \eta]+[\zeta, A(\eta)]\big)-[\zeta, \eta]\\
&=N_A(\zeta,\eta).
\end{split}
\end{equation}
So (cf.\ Example \ref{ia1} with $\beta=d^{\nabla} A,$ $\alpha=A$), $A$ is integrable iff 

\begin{equation*}
\begin{split}
\big(d^{\nabla} A \wedge (A \wedge A)\big)(\zeta,\eta)-d^{\nabla} A(\zeta,\eta)&=d^{\nabla} A(A(\zeta), A(\eta))-d^{\nabla} A(\zeta, \eta)\\
&=0.
\end{split}
\end{equation*}
\end{proof}

An expression similar to formula \ref{fog} above was obtained in \cite{Wan} (cf.\ the proof of Proposition 2.4). However, there is a mistake in the calculation (or at least in the English version, which is available as a preprint). This error seems to have been corrected in \cite{Wano} (cf.\ the proof of Proposition 3.2). The main distinction between Lemma \ref{Q1}, and the calculations from \cite{Wano, Wan} is that in these articles, the Levi-Civita connection was chosen in advance, while the computations here are valid for any torsion-free connection. The rewriting of \ref{fog} using the algebraic language of the previous section does not appear in \cite{Wano, Wan}, and is essentially new. The rewriting part is key as it allows for the proper covariant differentiation of the integrability form, $d^{\nabla} A \wedge (A \wedge A)-d^{\nabla} A,$ which will lead to curvature obstruction equations and curvature scalars as explained in the sequel. 

Intuitively, since the integrability form lives in $\Omega^2 (M, T_M),$ the correct notion of covariant differentiation is supplied by the covariant exterior derivative. In a basic way, since the integrability form is a function of $A$ and $d^{\nabla} A,$ its first covariant exterior derivative will be a function of $A,$ $d^{\nabla} A,$ and $(d^{\nabla})^2 A=R^{\nabla} \wedge A.$ The same thinking makes it evident that the $k$-th covariant exterior derivative of the integrability form will depend on \[A, d^{\nabla} A, (d^{\nabla})^2 A=R^{\nabla} \wedge A,\dots, (d^{\nabla})^{k+1} A.\] A formula for the last term is found in Lemma \ref{Q3}. The rest of this section will cover each step needed to derive the curvature obstruction equations.

Observe that 
\begin{lemma}{(Lemma 2, \cite{CSACPS})}\label{LRT}
If $\alpha \in \Omega^k(M,T_M),$ $\beta \in  \Omega^l(M,T_M),$ and if $\nabla'$ denotes the connection on $\bigwedge^{2} {T_M}$ that is induced by the torsion-free connection $\nabla$ on $T_M,$ then
\[d^{\nabla'}(\alpha \wedge \beta)=d^{\nabla} \alpha \wedge \beta+(-1)^k \alpha \wedge d^{\nabla} \beta.\]
\end{lemma}

Higher derivatives of the integrability form $d^{\nabla} A \wedge (A \wedge A)-d^{\nabla} A \in \Omega^2(M,T_M)$ will require a generalization of the above formula. But the first such derivative can already be computed, and it is included here as it serves as a guide for the more complicated cases.

\begin{remark}\label{excase}
Let $\alpha \in \Omega^1 (M, T_M),$ and put $I^{\nabla}_{\alpha}:=d^{\nabla} \alpha \wedge (\alpha \wedge \alpha)-d^{\nabla} \alpha.$ Then, \[d^{\nabla} I^{\nabla}_{\alpha}=(R^{\nabla} \wedge \alpha) \wedge (\alpha \wedge \alpha)-R^{\nabla} \wedge \alpha+2d^{\nabla} \alpha \wedge (d^{\nabla} \alpha \wedge \alpha).\] To see this, let $\nabla'$ be the induced connection on $\bigwedge^{2} {T_M}.$ By Lemma \ref{LRT} and the commuting rules for polyvector differential forms,
\begin{equation*}
\begin{split}
d^{\nabla} I^{\nabla}_{\alpha}&=\big((d^{\nabla})^2 \alpha \big) \wedge (\alpha \wedge \alpha)+d^{\nabla} \alpha \wedge \big(d^{\nabla'}(\alpha \wedge \alpha)\big)-(d^{\nabla})^2 \alpha \\
&=\big((d^{\nabla})^2 \alpha \big) \wedge (\alpha \wedge \alpha)-(d^{\nabla})^2 \alpha+d^{\nabla} \alpha \wedge \big(d^{\nabla} \alpha \wedge \alpha -\alpha \wedge d^{\nabla} \alpha\big)\\
&=(R^{\nabla} \wedge \alpha) \wedge (\alpha \wedge \alpha)-R^{\nabla} \wedge \alpha+2d^{\nabla} \alpha \wedge (d^{\nabla} \alpha \wedge \alpha).
 \end{split}
\end{equation*}
\end{remark}

The generalized formula for the $j$-th covariant exterior derivative w.r.t.\ the induced on $\bigwedge^{2} {T_M}$ connection $\nabla'$ will certainly be of the form

\[(d^{\nabla'})^j (A \wedge A)=\sum_{l=0}^j (-1)^l c_l^j \big((d^{\nabla})^{j-l} A\big) \wedge \big( (d^{\nabla})^l A\big),\] for some $c_l^j \in \mathbb{R}.$ From a direct computation of the first $10$ such derivatives emerges a pattern of generation of the unknown coefficients:

\[d^{\nabla'} (A \wedge A)=d^{\nabla} \wedge A - A \wedge d^{\nabla} A \implies c_0^1=c_1^1=1,\]

\[(d^{\nabla'})^2(A\wedge A)=\big((d^{\nabla})^2 A\big)\wedge A+d^{\nabla} A \wedge d^{\nabla} A-d^{\nabla} A \wedge d^{\nabla} A+A \wedge \big((d^{\nabla})^2 A\big) \implies c_0^2=1, c_1^2=0, c_2^2=1,\]

\[(d^{\nabla'})^3(A\wedge A)=\big((d^{\nabla})^3 A\big)\wedge A-\big((d^{\nabla})^2 A\big)\wedge d^{\nabla} A+d^{\nabla} A \wedge \big((d^{\nabla})^2 A\big)-A \wedge \big((d^{\nabla})^3 A\big)\] \[\implies c_0^3=c_1^3=c_2^3=c_3^3=1,\]

\[(d^{\nabla'})^4(A\wedge A)=\big((d^{\nabla})^4 A\big)\wedge A+2\big((d^{\nabla})^2 A\big)\wedge \big((d^{\nabla})^2 A\big)+A \wedge \big((d^{\nabla})^4 A\big)\] \[\implies c_0^4=1, c_1^4=0, c_2^4=2, c_3^4=0, c_4^4=1,\]

\begin{equation*}
\begin{split}
(d^{\nabla'})^5(A\wedge A)&=\big((d^{\nabla})^5 A\big)\wedge A-\big((d^{\nabla})^4 A\big)\wedge d^{\nabla} A+2\big((d^{\nabla})^3 A\big)\wedge \big((d^{\nabla})^2 A\big)-2\big((d^{\nabla})^2 A\big)\wedge \big((d^{\nabla})^3 A\big)+\\
&d^{\nabla} A \wedge \big((d^{\nabla})^4 A\big)-A \wedge \big((d^{\nabla})^5 A\big)
\end{split}
\end{equation*}
\[\implies c_0^5=1, c_1^5=1, c_2^5=2, c_3^5=2, c_4^5=1, c_5^5=1,\]

\begin{equation*}
\begin{split}
(d^{\nabla'})^6(A\wedge A)&=\big((d^{\nabla})^6 A\big)\wedge A+3\big((d^{\nabla})^4 A\big)\wedge \big((d^{\nabla})^2 A\big)+3\big((d^{\nabla})^2 A\big)\wedge \big((d^{\nabla})^4 A\big)+\\
&A \wedge \big((d^{\nabla})^6 A\big)
\end{split}
\end{equation*}
\[\implies c_0^6=1, c_1^6=0, c_2^6=3, c_3^6=0, c_4^6=3, c_5^6=0, c_6^6=1,\]

\begin{equation*}
\begin{split}
(d^{\nabla'})^7(A\wedge A)&=\big((d^{\nabla})^7 A\big)\wedge A-\big((d^{\nabla})^6 A\big)\wedge d^{\nabla} A+3\big((d^{\nabla})^5 A\big)\wedge \big((d^{\nabla})^2 A\big)+\\
&-3\big((d^{\nabla})^4 A\big)\wedge \big((d^{\nabla})^3 A\big)+3\big((d^{\nabla})^3 A\big)\wedge \big((d^{\nabla})^4 A\big)-3\big((d^{\nabla})^2 A\big)\wedge \big((d^{\nabla})^5 A\big)+\\
&d^{\nabla} A \wedge \big((d^{\nabla})^6 A\big)-A \wedge \big((d^{\nabla})^7 A\big)
\end{split}
\end{equation*}
\[\implies c_0^7=1, c_1^7=1, c_2^7=3, c_3^7=3, c_4^7=3, c_5^7=3, c_6^7=1,c_7^7=1,\]

\begin{equation*}
\begin{split}
(d^{\nabla'})^8(A\wedge A)&=\big((d^{\nabla})^8 A\big)\wedge A+4\big((d^{\nabla})^6 A\big)\wedge \big((d^{\nabla})^2 A\big)+6\big((d^{\nabla})^4 A\big)\wedge \big((d^{\nabla})^4 A\big)+\\
&4\big((d^{\nabla})^2 A\big)\wedge \big((d^{\nabla})^6 A\big)+A \wedge \big((d^{\nabla})^8 A\big)
\end{split}
\end{equation*}
\[\implies c_0^8=1, c_1^8=0, c_2^8=4, c_3^8=0, c_4^8=6, c_5^8=0, c_6^8=4, c_7^8=0,c_8^8=1,\]

\begin{equation*}
\begin{split}
(d^{\nabla'})^9(A\wedge A)&=\big((d^{\nabla})^9 A\big)\wedge A-\big((d^{\nabla})^8 A\big)\wedge d^{\nabla} A+4\big((d^{\nabla})^7 A\big)\wedge \big((d^{\nabla})^2 A\big)\\
&-4\big((d^{\nabla})^6 A\big)\wedge \big((d^{\nabla})^3 A\big)+6\big((d^{\nabla})^5 A\big)\wedge \big((d^{\nabla})^4 A\big)-6\big((d^{\nabla})^4 A\big)\wedge \big((d^{\nabla})^5 A\big)+\\
&4\big((d^{\nabla})^3 A\big)\wedge \big((d^{\nabla})^6 A\big)-4\big((d^{\nabla})^2 A\big)\wedge \big((d^{\nabla})^7 A\big)+d^{\nabla} A \wedge \big((d^{\nabla})^8 A\big)\\
&-A \wedge \big((d^{\nabla})^9 A\big)
\end{split}
\end{equation*}
\[\implies c_0^9=1, c_1^9=1, c_2^9=4, c_3^9=4, c_4^9=6, c_5^9=6, c_6^9=4, c_7^9=4,c_8^9=1, c_9^9=1,\]

\begin{equation*}
\begin{split}
(d^{\nabla'})^{10}(A\wedge A)&=\big((d^{\nabla})^{10} A\big)\wedge A+5\big((d^{\nabla})^8 A\big)\wedge \big((d^{\nabla})^2 A\big)+10\big((d^{\nabla})^6 A\big)\wedge \big((d^{\nabla})^4 A\big)+\\
&10\big((d^{\nabla})^4 A\big)\wedge \big((d^{\nabla})^6 A\big)+5\big((d^{\nabla})^2 A\big)\wedge \big((d^{\nabla})^8 A\big)+A \wedge \big((d^{\nabla})^{10} A\big)
\end{split}
\end{equation*}
\[\implies c_0^{10}=1, c_1^{10}=0, c_2^{10}=5, c_3^{10}=0, c_4^{10}=10, c_5^{10}=0, c_6^{10}=10, c_7^{10}=0,c_8^{10}=5, c_9^{10}=0, c_{10}^{10}=1.\]

And one trivially has that \[(d^{\nabla'})^0(A\wedge A)=A \wedge A \implies c_0^0=1.\]

\begin{figure}
\caption{$c_l^j$ coefficients}
\begin{tabular}{>{$j=}l<{$\hspace{12pt}}*{21}{c}}
0 &&&&&&&&&&&1&&&&&&&&&&\\
1 &&&&&&&&&&1&&1&&&&&&&&&\\
2 &&&&&&&&&1&&0&&1&&&&&&&&\\
3 &&&&&&&&1&&1&&1&&1&&&&&&&\\
4 &&&&&&&1&&0&&2&&0&&1&&&&&&\\
5 &&&&&&1&&1&&2&&2&&1&&1&&&&&\\
6 &&&&&1&&0&&3&&0&&3&&0&&1&&&&\\
7 &&&&1&&1&&3&&3&&3&&3&&1&&1&&&\\
8 &&&1&&0&&4&&0&&6&&0&&4&&0&&1&&\\
9 &&1&&1&&4&&4&&6&&6&&4&&4&&1&&1&\\
10 &1&&0&&5&&0&&10&&0&&10&&0&&5&&0&&1
\end{tabular}
\label{faust}
\end{figure}

So one way to describe this pattern of coefficients is by saying that if $j$ is even, then \[c_l^j=\begin{cases} {j/2 \choose l/2} & l \mbox{ is even}\\
0 & l \mbox{ is odd},
 \end{cases}\]
 while if $j$ is odd, the corresponding coefficients are recursively given by $c_l^j=c_{l-1}^{j-1}+c_l^{j-1}.$ And this numerical data can be visually arranged in a pyramidal fashion akin to a hybrid of the Pascal and Sierpinski triangles as illustrated in Figure \ref{faust}. 
\begin{lemma}{(A general Leibniz rule)}\label{AGLR}
For $j\geq 1,$

\[(d^{\nabla'})^j(A \wedge A)=\sum_{l=0}^j \frac{1}{2^\frac{{(-1)^j}+1}{2}}\Big(\frac{(-1)^j+1}{2}+(-1)^l \Big){\lfloor{\frac{j}{2}\rfloor} \choose \lfloor{\frac{l}{2}}\rfloor} \big((d^{\nabla})^{j-l} A\big) \wedge \big((d^{\nabla})^l A\big),\] where  $\lfloor{ \cdot }\rfloor$ is the floor function. More explicitly, for $m \geq 0,$

\begin{equation}\label{vv1}
\begin{split}
(d^{\nabla'})^{2m} (A \wedge A)=\sum_{l=0}^m {m \choose l} \big((d^{\nabla})^{2m-2l} A\big) \wedge \big((d^{\nabla})^{2l} A\big),
\end{split}
\end{equation}
and
\begin{equation}\label{vv2}
\begin{split}
(d^{\nabla'})^{2m+1} (A \wedge A)=\sum_{l=0}^m {m \choose l} \Big[\big((d^{\nabla})^{2m+1-2l} A\big) \wedge \big((d^{\nabla})^{2l} A\big)-\big((d^{\nabla})^{2m-2l} A\big) \wedge \big((d^{\nabla})^{2l+1} A\big)\Big].
\end{split}
\end{equation}
\end{lemma}

\begin{proof}
The proof is by induction. But first, observe that assuming the validity of the main equation claimed, formulas \ref{vv1} and \ref{vv2} follow at once by taking $j=2m,$ and $j=2m+1.$

When $j=1,$ 
\begin{equation*}
\begin{split}
d^{\nabla'} (A\wedge A)&=\sum_{l=0}^1 (-1)^l {0 \choose \lfloor \frac{l}{2} \rfloor} \big((d^{\nabla})^{1-l} A\big) \wedge \big((d^{\nabla})^l A\big)\\
&=d^{\nabla} A \wedge A - A \wedge d^{\nabla} A,
\end{split}
\end{equation*}
consistent with Lemma \ref{LRT}. 

For the induction hypothesis, assume that formulas \ref{vv1} and \ref{vv2} are true for $m.$ Notice that
\begin{equation*}
\begin{split}
(d^{\nabla'})^{2m+1} (A \wedge A)&=(d^{\nabla'})\big( (d^{\nabla'})^{2m} (A \wedge A)\big)\\
&=(d^{\nabla'}) \Big[\sum_{l=0}^m {m \choose l} \big((d^{\nabla})^{2m-2l} A\big) \wedge \big((d^{\nabla})^{2l} A\big) \Big]\\
&=\sum_{l=0}^m {m \choose l} \Big[\big((d^{\nabla})^{2m+1-2l} A\big) \wedge \big((d^{\nabla})^{2l} A\big)-\big((d^{\nabla})^{2m-2l} A\big) \wedge \big((d^{\nabla})^{2l+1} A\big)\Big],
\end{split}
\end{equation*}

where the last line follows from Lemma \ref{LRT}, and 

\begin{equation*}
\begin{split}
(d^{\nabla'})^{2m+2}&=d^{\nabla'} \big((d^{\nabla'})^{2m+1} (A \wedge A) \big)\\
&=d^{\nabla'} \Big\{\sum_{l=0}^m {m \choose l} \Big[\big((d^{\nabla})^{2m+1-2l} A\big) \wedge \big((d^{\nabla})^{2l} A\big)-\big((d^{\nabla})^{2m-2l} A\big) \wedge \big((d^{\nabla})^{2l+1} A\big)\Big] \Big\}\\
&=\sum_{l=0}^m {m \choose l} \Big[\big((d^{\nabla})^{2m+2-2l} A\big)\wedge \big((d^{\nabla})^{2l} A\big)+\big((d^{\nabla})^{2m+1-2l} A\big)\wedge \big((d^{\nabla})^{2l+1} A\big)\\
&-\big((d^{\nabla})^{2m-2l+1} A\big)\wedge \big((d^{\nabla})^{2l+1} A\big)+\big((d^{\nabla})^{2m-2l} A\big)\wedge \big((d^{\nabla})^{2l+2} A\big) \Big]\\
&=\sum_{l=0}^m \Big[\big((d^{\nabla})^{2m+2-2l} A\big)\wedge \big((d^{\nabla})^{2l} A\big)+\big((d^{\nabla})^{2m-2l} A\big)\wedge \big((d^{\nabla})^{2l+2} A\big) \Big]\\
&=\sum_{l=0}^m {m \choose l} \big((d^{\nabla})^{2m+2-2l} A\big)\wedge \big((d^{\nabla})^{2l} A\big)+\sum_{l=1}^{m+1} {m \choose l-1} \big((d^{\nabla})^{2m-2l+2} A\big)\wedge \big((d^{\nabla})^{2l} A\big)\\
&=\sum_{l=0}^{m+1} \Big[{m \choose l}+{m \choose l-1} \Big]\big((d^{\nabla})^{2m+2-2l} A\big)\wedge \big((d^{\nabla})^{2l} A\big)\\
&=\sum_{l=0}^{m+1} {m+1 \choose l} \big((d^{\nabla})^{2m+2-2l} A\big)\wedge \big((d^{\nabla})^{2l} A\big),
\end{split}
\end{equation*}
so the proof is complete.
\end{proof}

The $k$-th exterior covariant derivative of the integrability form is geometrically meaningful, and this will start becoming apparent in Lemma \ref{Q3}.

\begin{lemma}\label{Q3}
For any $A \in \Omega^1 (M,T_M),$ 

\[(d^{\nabla})^{2m} A =(R^{\nabla})^m \wedge A,\mbox{ and }(d^{\nabla})^{2m+1} A =(R^{\nabla})^m \wedge d^{\nabla} A,\] where
\[\big((R^{\nabla})^m \wedge A\big)(X_1,\dots,X_{2m+1})=\frac{1}{(2m)!} \sum_{\sigma \in S_{2m+1}} sign(\sigma)(R^{\nabla})^m (X_{\sigma(1)},\dots,X_{\sigma(2m)})\big(A(X_{\sigma(2m+1)})\big),\] and 
\begin{equation*}
\begin{split}
\big((R^{\nabla})^m \wedge d^{\nabla} A\big)(X_1,\dots,X_{2m+2})&=\frac{1}{2(2m)!} \sum_{\sigma \in S_{2m+2}} sign(\sigma)(R^{\nabla})^m (X_{\sigma(1)},\dots,X_{\sigma(2m)})\big(d^{\nabla} A (X_{\sigma(2m+1)},\\
& X_{\sigma(2m+2)})\big).
\end{split}
\end{equation*}
\end{lemma}

\begin{proof}
By induction: the even base case $(d^{\nabla})^2 A=R^{\nabla} \wedge A$ can be verified using the definition of $d^{\nabla}$; indeed, 
\begin{equation*}
\begin{split}
\big((d^{\nabla})^2 A\big)(X,Y,Z)&=R^{\nabla}(X,Y)(A(Z))+R^{\nabla}(Y,Z)(A(X))+R^{\nabla}(Z,X)(A(Y))\\
&=\big(R^{\nabla} \wedge A\big)(X,Y,Z).
\end{split}
\end{equation*}

Let $\nabla'$ be the induced connection on $\End_{\mathbb{R}}(T_M),$ and recall that $d^{\nabla'} R^{\nabla}=0.$ By the above,
\begin{equation*}
\begin{split}
(d^{\nabla})^3 A&=d^{\nabla} (R^{\nabla} \wedge A)\\
&=R^{\nabla} \wedge d^{\nabla} A,
\end{split}
\end{equation*}
establishing the odd case. 

Now, if $(d^{\nabla})^{2k} A =(R^{\nabla})^k \wedge A,$ then

\begin{equation*}
\begin{split}
(d^{\nabla})^{2k+2} A&=(d^{\nabla})^2 \big((R^{\nabla})^k \wedge A\big)\\
&=(R^{\nabla})^k \wedge \big((d^{\nabla})^2 A\big) \\
&=(R^{\nabla})^{k+1} \wedge A.
\end{split}
\end{equation*}
For the odd case, assuming that $(d^{\nabla})^{2k+1} A =(R^{\nabla})^k \wedge d^{\nabla} A,$ the above even case implies that 
\begin{equation*}
\begin{split}
(d^{\nabla})^{2k+3} A&=d^{\nabla} \big((d^{\nabla})^{2k+2} A\big)\\
&=d^{\nabla} \big((R^{\nabla})^{k+1} \wedge A \big) \\
&=(R^{\nabla})^{k+1} \wedge d^{\nabla} A
\end{split}
\end{equation*}
as desired.
\end{proof}

\begin{lemma}\label{wocurv}
For any $m \geq 0,$ 

\begin{equation*}
\begin{split}
(d^{\nabla})^{2m} \big(d^{\nabla} A \wedge (A \wedge A)\big)&=\sum_{j=0}^m \sum_{l=0}^j \Big[{2m \choose 2j} {j \choose l} \big((d^{\nabla})^{2(m-j)+1} A \big)\wedge \Big( \big((d^{\nabla})^{2(j-l)} A\big) \wedge \big((d^{\nabla})^{2l} A\big)\Big)+\\
&{2m \choose 2j+1} {j \choose l} \big((d^{\nabla})^{2(m-j)} A \big)\wedge \Big( \big((d^{\nabla})^{2(j-l)+1} A\big) \wedge \big((d^{\nabla})^{2l} A\big)\Big)\\
&-{2m \choose 2j+1} {j \choose l} \big((d^{\nabla})^{2(m-j)} A \big)\wedge \Big( \big((d^{\nabla})^{2(j-l)} A\big) \wedge \big((d^{\nabla})^{2l+1} A\big)\Big)\Big],
\end{split}
\end{equation*}
and 
\begin{equation*}
\begin{split}
(d^{\nabla})^{2m+1} \big(d^{\nabla} A \wedge (A \wedge A)\big)&=\sum_{j=0}^m \sum_{l=0}^j \Big[{2m+1 \choose 2j} {j \choose l} \big((d^{\nabla})^{2(m-j+1)} A \big)\wedge \Big( \big((d^{\nabla})^{2(j-l)} A\big) \wedge \big((d^{\nabla})^{2l} A\big)\Big)+\\
&{2m+1 \choose 2j+1} {j \choose l} \big((d^{\nabla})^{2(m-j)+1} A \big)\wedge \Big( \big((d^{\nabla})^{2(j-l)+1} A\big) \wedge \big((d^{\nabla})^{2l} A\big)\Big)\\
&-{2m+1 \choose 2j+1} {j \choose l} \big((d^{\nabla})^{2(m-j)+1} A \big)\wedge \Big( \big((d^{\nabla})^{2(j-l)} A\big) \wedge \big((d^{\nabla})^{2l+1} A\big)\Big)\Big],
\end{split}
\end{equation*}
\end{lemma}

\begin{proof}
Observe that \[(d^{\nabla})^k \big(d^{\nabla} A \wedge (A \wedge A)\big)=\sum_{j=0}^k {k \choose j} \big((d^{\nabla})^{k+1-j} A \big) \wedge \big((d^{\nabla'})^j (A \wedge A)\big).\] The calculations follow from organizing the terms in each sum according to parity, and an application of Lemma \ref{AGLR}.
\end{proof}

The integrability of an almost-complex structure constrains the curvature of any torsion-free connection $\nabla$ on $T_M$ in the following sense. 

\begin{proposition}\label{PPP1}
Suppose that $A \in AC(M)$ is integrable. Then, for $m \geq 0,$ 
\begin{equation}\label{equ1}
\begin{split}
&\sum_{j=0}^m \sum_{l=0}^j \Big[{2m \choose 2j} {j \choose l} \big((R^{\nabla})^{m-j}\wedge d^{\nabla} A \big)\wedge \Big( \big((R^{\nabla})^{j-l} \wedge A\big) \wedge \big((R^{\nabla})^l \wedge A\big)\Big)+\\
&{2m \choose 2j+1} {j \choose l} \big((R^{\nabla})^{m-j} \wedge A \big)\wedge \Big( \big((R^{\nabla})^{j-l} \wedge d^{\nabla}A\big) \wedge \big((R^{\nabla})^l \wedge A\big)\Big)\\
&-{2m \choose 2j+1} {j \choose l} \big((R^{\nabla})^{m-j} \wedge A \big)\wedge \Big( \big((R^{\nabla})^{j-l} \wedge A\big) \wedge \big((R^{\nabla})^{l} \wedge d^{\nabla} A\big)\Big)\Big]-(R^{\nabla})^m \wedge d^{\nabla} A=0,
\end{split}
\end{equation}
and
\begin{equation}\label{equ2}
\begin{split}
&\sum_{j=0}^m \sum_{l=0}^j \Big[{2m+1 \choose 2j} {j \choose l} \big((R^{\nabla})^{m-j+1}\wedge A \big)\wedge \Big( \big((R^{\nabla})^{j-l} \wedge A\big) \wedge \big((R^{\nabla})^l \wedge A\big)\Big)+\\
&{2m+1 \choose 2j+1} {j \choose l} \big((R^{\nabla})^{m-j} \wedge d^{\nabla} A \big)\wedge \Big( \big((R^{\nabla})^{j-l} \wedge d^{\nabla}A\big) \wedge \big((R^{\nabla})^l \wedge A\big)\Big)\\
&-{2m+1 \choose 2j+1} {j \choose l} \big((R^{\nabla})^{m-j} \wedge d^{\nabla} A \big)\wedge \Big( \big((R^{\nabla})^{j-l} \wedge A\big) \wedge \big((R^{\nabla})^{l} \wedge d^{\nabla} A\big)\Big)\Big]-(R^{\nabla})^{m+1} \wedge A=0.
\end{split}
\end{equation}
\end{proposition}

\begin{proof}
By Lemma \ref{Q1}, if $A$ is integrable, then any covariant exterior derivative of the form $d^{\nabla} A \wedge (A \wedge A)-d^{\nabla} A$ must vanish. The left hand side of equation \ref{equ1} is \[(d^{\nabla})^{2m} \big(d^{\nabla} A \wedge (A \wedge A)-d^{\nabla} A\big)\] but rewritten using Lemmas \ref{Q3} and \ref{wocurv}. The same reasoning works for verifying equation \ref{equ2} whose left hand side is instead  \[(d^{\nabla})^{2m+1} \big(d^{\nabla} A \wedge (A \wedge A)-d^{\nabla} A\big).\] 
\end{proof}

\begin{remark}\label{REM1}    
The first 4 curvature obstruction equations are

\begin{equation*}
(R^{\nabla} \wedge A)\wedge (A \wedge A)+2d^{\nabla} A \wedge (d^{\nabla} A \wedge A)-R^{\nabla} \wedge A=0,
\end{equation*} 

\begin{equation*}
\begin{split}
&(R^{\nabla} \wedge d^{\nabla} A)\wedge (A \wedge A)+4(R^{\nabla} \wedge A)\wedge (d^{\nabla} A \wedge A)+2d^{\nabla} A\wedge \big((R^{\nabla} \wedge A)\wedge A\big)-R^{\nabla} \wedge d^{\nabla} A=0,
\end{split}
\end{equation*}

\begin{equation*}
\begin{split}
&\big((R^{\nabla})^2 \wedge A\big)\wedge (A \wedge A)+6(R^{\nabla} \wedge d^{\nabla} A)\wedge (d^{\nabla} A \wedge A)+6(R^{\nabla} \wedge A)\wedge \big((R^{\nabla} \wedge A)\wedge A\big)+\\
&2 d^{\nabla} A \wedge \big((R^{\nabla} \wedge d^{\nabla} A)\wedge  A\big)-2d^{\nabla} A \wedge \big((R^{\nabla} \wedge A)\wedge d^{\nabla} A\big)-(R^{\nabla})^2 \wedge A=0,
\end{split}
\end{equation*}
and
\begin{equation*}
\begin{split}
&\big((R^{\nabla})^2 \wedge d^{\nabla} A\big)\wedge (A \wedge A)+8\big((R^{\nabla})^2 \wedge A\big)\wedge (d^{\nabla} A \wedge A)+12(R^{\nabla} \wedge d^{\nabla} A)\wedge \big((R^{\nabla} \wedge A)\wedge A\big)+\\
&8(R^{\nabla} \wedge A) \wedge \big((R^{\nabla} \wedge d^{\nabla} A)\wedge  A\big)-8(R^{\nabla} \wedge A)\wedge \big((R^{\nabla} \wedge A)\wedge d^{\nabla} A\big)+\\
&2d^{\nabla} A\wedge \Big(\big((R^{\nabla})^2 \wedge A\big)\wedge A \Big)+2d^{\nabla} A\wedge \big((R^{\nabla} \wedge A)\wedge (R^{\nabla} \wedge A)\big)-(R^{\nabla})^2 \wedge d^{\nabla} A=0.
\end{split}
\end{equation*}
\end{remark}

\begin{proof}
The computation of $(d^{\nabla'})^j (A. \wedge A),$ $j=1,\dots,4,$ further simplifies by means of the appropriate commuting laws:

\[d^{\nabla'} (A \wedge A)=2d^{\nabla} A\wedge A,\]

\[(d^{\nabla'})^2(A\wedge A)=2\big((d^{\nabla})^2 A\big)\wedge A,\]

\[(d^{\nabla'})^3(A\wedge A)=2\big((d^{\nabla})^3 A\big)\wedge A-2\big((d^{\nabla})^2 A\big)\wedge d^{\nabla} A,\]

\[(d^{\nabla'})^4(A\wedge A)=2\big((d^{\nabla})^4 A\big)\wedge A+2\big((d^{\nabla})^2 A\big)\wedge \big((d^{\nabla})^2 A\big).\]

These computations come into play when taking covariant exterior derivatives of the integrability form $d^{\nabla} A \wedge (A \wedge A)-d^{\nabla} A,$ and thereby lead to curvature obstruction equations through Lemma \ref{Q3} as shown below.

The first such derivative is

\begin{equation*}
\begin{split}
d^{\nabla} \big(d^{\nabla} A \wedge (A \wedge A)-d^{\nabla} A\big)&=\sum_{j=0}^1 {1 \choose j} \big((d^{\nabla})^{2-j} A\big) \wedge \big((d^{\nabla'})^j (A\wedge A\big)-\big((d^{\nabla})^2 A\big)\\
&=\big((d^{\nabla})^{2} A\big) \wedge (A \wedge A)+2d^{\nabla} A \wedge (d^{\nabla} A \wedge A)-\big((d^{\nabla})^2 A\big)\\
&=(R^{\nabla} \wedge A)\wedge (A \wedge A)+2d^{\nabla} A \wedge (d^{\nabla} A \wedge A)-R^{\nabla} \wedge A,
\end{split}
\end{equation*}
and the second one is
\begin{equation*}
\begin{split}
(d^{\nabla})^2 \big(d^{\nabla} A \wedge (A \wedge A)-d^{\nabla} A\big)&=\big((d^{\nabla})^{3} A\big) \wedge (A \wedge A)+4\big((d^{\nabla})^2 A\big) \wedge (d^{\nabla} A \wedge A)+2d^{\nabla} A\wedge \Big(\big((d^{\nabla})^2 A\big) \wedge A\Big)\\
&-\big((d^{\nabla})^3 A\big)\\
&=(R^{\nabla} \wedge d^{\nabla} A)\wedge (A \wedge A)+4(R^{\nabla} \wedge A)\wedge (d^{\nabla} A \wedge A)+\\
&2d^{\nabla} A\wedge \big((R^{\nabla} \wedge A)\wedge A\big)-R^{\nabla} \wedge d^{\nabla} A.
\end{split}
\end{equation*}
Hence, the first $2$ curvature obstruction equations are \[(R^{\nabla} \wedge A)\wedge (A \wedge A)+2d^{\nabla} A \wedge (d^{\nabla} A \wedge A)-R^{\nabla} \wedge A=0,\] and 
\begin{equation*}
\begin{split}
&(R^{\nabla} \wedge d^{\nabla} A)\wedge (A \wedge A)+4(R^{\nabla} \wedge A)\wedge (d^{\nabla} A \wedge A)+2d^{\nabla} A\wedge \big((R^{\nabla} \wedge A)\wedge A\big)-R^{\nabla} \wedge d^{\nabla} A=0.
\end{split}
\end{equation*}
The order $k,$ for $3 \leq k \leq 4,$ derivatives of the integrability form can be computed the same way, leading to the remaining obstruction equations.
\end{proof}

By choosing to work with the Levi-Civita connection, one may distill the Riemann, Ricci, sectional, and scalar curvature, versions of these obstruction equations, among others. It could be interesting to investigate those equations and their applications to the complex structure existence (CSE) problem. The contracted obstruction equations could be helpful for studying the CSE problem for almost-complex manifolds admitting metrics with bounded curvature. 

\section{Constant curvature}\label{csc}
Attention will now be placed on a distinguished family of complex structures, those $A \in C(M)$ such that $d^{\nabla} A=0,$ for a specific symmetric connection $\nabla$ on $T_M.$ Similar structures have been studied in, for example, \cite{Free}. These special complex structures hold space to reflect on the relationship between the Riemannian and (almost-)complex geometries of $M.$ 

Recall that if $\nabla$ is the Levi-Civita connection of a Riemannian metric $g$ on $M,$ then the Riemann curvature tensor $Rm$ is given as \[Rm(X,Y,Z,W)=g\Big(R^{\nabla} (X,Y)Z, W\Big).\] For any $p \in M,$ $2$-plane $\Pi_p \subset T_{M,p},$ and basis $(X,Y)$ of $\Pi_p,$ \[K(\Pi_p)=K(X,Y):=\frac{Rm(X,Y,X,Y)}{g(X,X)g(Y,Y)-g(X,Y)^2}\] is the sectional curvature at $p$ with respect to $\Pi_p.$ Then, $(M,g)$ is said to have constant sectional curvature $c_0$ if $K(\Pi_p)=c_0$ for all $p \in M,$ and all $2$-planes $\Pi_p \subset T_{M,p}.$

The Kulkarni-Nomizu product of two $(0,2)$ symmetric tensors $h$ and $k$ is defined as \[h \owedge k (X,Y,Z,W):=h(X,Z)k(Y,W)+h(Y,W)k(X,Z)-h(X,W)k(Y,Z)-h(Y,Z)k(X,W).\] The Riemann curvature tensor of a metric $g$ of constant sectional curvature $c_0$ is of the form \[Rm=-\frac{c_0}{2} g \owedge g.\] 

\begin{lemma}\label{LT1}
Let $X, Y, Z, W \in \mathfrak{X}(M),$ and $A \in \Omega^1 (M, T_M)$ be such that $d^{\nabla} A=0.$ Then,
\begin{enumerate}
\item  $R^{\nabla}(X,Y)(A(Z))+R^{\nabla}(Y,Z)(A(X))+R^{\nabla}(Z,X)(A(Y))=0,$ and
\item in the special case when $\nabla$ is the Levi-Civita connection of a Riemannian metric on $M,$  \[Rm(X,Y,A(Z),W)+Rm(Y,Z,A(X),W)+Rm(Z,X,A(Y),W)=0.\]
\end{enumerate}
\end{lemma}

\begin{proof}
The first claim is immediate from $0=(d^{\nabla})^2 A=R^{\nabla} \wedge A,$ and the second claim, from the definition of $Rm$ and 1.
\end{proof}

\begin{theorem}\label{TV1}
Let $(M,g)$ be a Riemannian manifold of real dimension at least $4.$ If $g$ has non-zero constant sectional curvature, then $M$ does not admit a complex structure $A \in C(M)$ satisfying $d^{\nabla} A=0,$ where $\nabla$ is the Levi-Civita connection.
\end{theorem}

\begin{proof}
Suppose that $g$ has constant sectional curvature $c_0 \neq 0.$ To the contrary, assume that there exists $A \in AC(M)$ that satisfies $d^{\nabla} A=0,$ for the Levi-Civita connection $\nabla$ on $T_M.$ Note that $A$ is then automatically integrable (cf.\ Lemma \ref{Q1}).

Let $p \in M,$ and $(E_i)_{i=1}^n$ be an orthonormal basis of $T_{M,p}.$ Since $Rm=-\frac{c_0}{2} g \owedge g,$  \[R_{ijkl}:=Rm(E_i,E_j,E_k,E_l)=c_0(\delta_{il}\delta_{jk}-\delta_{ik}\delta_{jl}),\] and so, if $i\neq j,$ $R_{ijji}=-R_{ijij}=c_0,$ all other components of $Rm$ being equal to zero. Lemma \ref{LT1}, 2.\ with $X=W$ says that
\begin{equation}\label{eq3}
Rm(X,Y,A(Z),X)+Rm(Y,Z,A(X),X)+Rm(Z,X,A(Y),X)=0.
\end{equation}
At $p,$ with $X_p=E_i,$ $Y_p=E_j,$ $Z_p=E_k$ so that $A(E_a)=\sum_{b=1}^n A^b_{a} E_b,$ equation \ref{eq3} becomes 

\begin{equation*}
\sum_{b=1}^n \Big[A^b_{k} R_{ijbi}+A^b_{i}R_{jkbi}+A^b_{j}R_{kibi} \Big]=0.
\end{equation*}

Now, impose the requirement that $i,j,k$ be distinct integers between $1$ and $n$ (this is why $n$ is assumed to be $4$ or larger in the hypothesis). Then, 

\begin{equation*}
\begin{split}
0&=\sum_{b=1}^n \Big[A^b_{k} R_{ijbi}+A^b_{i}R_{jkbi}+A^b_{j}R_{kibi} \Big]\\
&=A^j_{k} R_{ijji}+A^k_{j}R_{kiki}\\
&=c_0 (A^j_{k}-A^k_{j}).
\end{split}
\end{equation*}
This is a contradiction because an almost-complex structure $A \in AC(M)$ at a point $p$ cannot be symmetric -- the spectral theorem would imply that $A$ has real eigenvalues.
\end{proof}

Let \[AC(M)_g:=\{A\in AC(M)\mid g(A\zeta,A\eta)=g(\zeta,\eta), \forall \zeta,\eta \in \mathfrak{X}(M)\}\] be the space of almost-complex structures that are compatible with a Riemannian metric $g.$ Recall that a K{\"a}hler manifold is a triple $(M,A,g)$ of complex manifold $(M,A)$ with $A \in C(M) \cap AC(M)_g,$ where $g$ is a Riemannian metric for which the fundamental form $\Omega_A:=g(A\cdot,\cdot)$ is closed. Equivalently, $(M,A,g)$ with $A \in AC(M)_g$ (i.e.\ almost-hermitian) is a K{\"a}hler iff $\nabla A=0,$ where $\nabla$ is the Levi-Civita connection of $g.$ 

Certainly, for any torsion-free connection $\nabla$ on $T_M,$ and $A \in \Omega^1(M,T_M),$ $\nabla A=0$ implies that $d^{\nabla} A=0.$ The converse is true for almost-hermitian $(M,A,g),$ and the Levi-Civita connection. As a result, the almost-hermitian manifold $(M,A,g)$ is K{\"a}hler iff $d^{\nabla} A=0,$ where $\nabla$ is the Levi-Civita connection. In view of Theorem \ref{TV1},

\begin{remark}\label{R1}
The underlying Riemannian metric of a K{\"a}hler manifold $(M,A,g)$ of real dimension at least $4$ cannot have non-zero constant sectional curvature.
\end{remark}

\section{Functionals and numerical obstructions}\label{funct}
The material presented in the first part of this section is closely related to \cite{CSACPS}. Let $g$ be a Riemannian metric on $M,$ and $(x_i)_{i=1}^n$ be local coordinates, giving the associated local frame $\big(\frac{\partial}{\partial x_i}\big)_{i=1}^n$ of $T_M.$ Any element $\rho \in \Omega^k(M,T_M)$ can be locally expressed as $\rho=\sum_{i=1}^n p_i \otimes \frac{\partial}{\partial x_i}.$ There is a natural $L^2$-inner product on $\Omega^{\bullet}(M,T_M)$ that comes from a natural extension of the Hodge star operator to $\Omega^{\bullet}(M,T_M).$ Indeed, if \[\wedge_g:\Omega^k(M,T_M) \otimes \Omega^l(M,T_M) \to \Omega^{k+l}(M)\] is the bilinear map \[\alpha \wedge_g \beta=\sum_{i,j=1}^n a_i \wedge b_j g\big(\frac{\partial}{\partial x_i}, \frac{\partial}{\partial x_j}\big),\] and \[\langle \cdot, \cdot \rangle_g:\Omega^k(M,T_M) \otimes \Omega^k(M,T_M) \to \mathbb{R}\] is the pairing \[\langle \alpha, \alpha'\rangle_g=\sum_{i,j=1}^n \langle a_i, a'_j \rangle g\big(\frac{\partial}{\partial x_i}, \frac{\partial}{\partial x_j}\big),\] the induced Hodge star operator \[\star_g:\Omega^k(M,T_M) \to \Omega^{n-k}(M,T_M)\] is defined by the rule \[\alpha \wedge_g \star_g \alpha'=\langle \alpha, \alpha'\rangle_g \vol_g.\] Then, the $L^2$-product is \[\langle \langle \cdot, \cdot \rangle \rangle : \Omega^{\bullet} (M, T_M) \otimes \Omega^{\bullet} (M, T_M) \to \mathbb{R},\] \[\langle \langle A, B \rangle \rangle :=\sum_k \langle A_k, B_k \rangle_k,\] where \[\langle A_k, B_k \rangle_k:=\int_M A_k \wedge_g \star_g B_k,\] and where $A=\sum_k A_k,$ $B=\sum_k B_k,$ $A_k, B_k \in \Omega^k (M, T_M).$ Let $\delta^{\nabla}$ be the formal adjoint of $d^{\nabla}$ w.r.t.\ this inner product.

\begin{proposition}\label{PQ}
Let $M$ be a real manifold, and $\nabla$ be a symmetric connection on $T_M.$ There exist functionals, denoted here by $\mathcal{C}^{\nabla}[1]$ and $\mathcal{C}^{\nabla}[1, 3],$ whose sets of critical points, when projected onto the degree $1$ piece of $\Omega^{\bullet} (M, T_M),$ coincide with $\ker{\big(d^{\nabla}\big|_{\Omega^1(M,T_M)}\big)}.$
\end{proposition}

\begin{proof}
Let $p_k:\Omega^{\bullet} (M, T_M) \to \Omega^k(M,T_M)$ be the projection mapping $p_k(\gamma)=\gamma_k.$ Adjust the covariant exterior derivative associated with the chosen connection by setting \[d^{\nabla}[1]:=\sum_{i=0}^{n-1} (1-\delta_0^i)d^{\nabla}\big|_{\Omega^i(M,T_M)}.\] This new operator coincides with $d^{\nabla}$ on each degree $i$ piece of $\Omega^{\bullet} (M, T_M),$ except for $i=0,$ when it is zero. The formal adjoint, $\delta^{\nabla}[1],$ then agrees with $\delta^{\nabla}$ on $\Omega^i(M,T_M),$ for all $i \neq 1,$ and is zero for $i=1.$ 

Put \[\widetilde{\Omega}^1(M,T_M):=\{\gamma \in \Omega^{\bullet} (M, T_M) \mid \delta^{\nabla} \gamma_3=0\}.\] Consider the functional \[\mathcal{C}^{\nabla}[1]:\widetilde{\Omega}^1(M,T_M) \to \mathbb{R}, \quad \mathcal{C}^{\nabla}[1](\gamma):=\langle\langle d^{\nabla}[1] \gamma, \gamma\rangle \rangle.\] The first variation of this functional is \[\frac{d}{dt}\Big|_{t=0} \mathcal{C}^{\nabla}[1](\gamma+t\beta)=\langle\langle \beta, (\delta^{\nabla}[1]+d^{\nabla}[1])\gamma\rangle\rangle,\] so the set of critical points is \[\mathcal{S}=\{\gamma \in \widetilde{\Omega}^1(M,T_M) \mid \delta^{\nabla}[1]\gamma_{k+1}+d^{\nabla}[1]\gamma_{k-1}=0, \forall k \neq 2, \mbox{ and }d^{\nabla} \gamma_1=0\}.\] Thus, $p_1(\mathcal{S})=\ker{\big(d^{\nabla}\big|_{\Omega^1(M,T_M)}\big)}.$

Alternatively, (re-)define the covariant exterior derivative as \[d^{\nabla}[1,3]:=\sum_{i=0}^{n-1} (1-\delta_0^i-\delta_2^i)d^{\nabla}\big|_{\Omega^i(M,T_M)}.\] And now consider the functional \[\mathcal{C}^{\nabla}[1,3]:\Omega^{\bullet} (M, T_M) \to \mathbb{R}, \quad \mathcal{C}^{\nabla}[1, 3](\gamma):=\langle\langle d^{\nabla}[1, 3] \gamma, \gamma\rangle \rangle.\] The first variation is \[\frac{d}{dt}\Big|_{t=0} \mathcal{C}^{\nabla}[1, 3](\gamma+t\beta)=\langle\langle \beta, (\delta^{\nabla}[1, 3]+d^{\nabla}[1, 3])\gamma\rangle\rangle.\] The critical point set is \[\mathcal{S}'=\{\gamma \in \Omega^{\bullet} (M, T_M) \mid \delta^{\nabla}[1, 3]\gamma_{k+1}+d^{\nabla}[1, 3]\gamma_{k-1}=0, \forall k \neq 2, \mbox{ and }d^{\nabla} \gamma_1=0\},\] so $p_1(\mathcal{S}')=\ker{\big(d^{\nabla}\big|_{\Omega^1(M,T_M)}\big)}$ as well. 
\end{proof}

\begin{corollary}\label{CRO}
Let $\widetilde{\mathcal{C}}(M):=\{\gamma \in \Omega^{\bullet} (M, T_M) \mid \gamma_1 \in C(M)\}.$ The special complex structures $A \in C(M)$ with the added property that $d^{\nabla} A=0$ (cf.\ Theorem \ref{TV1}) can be variationally realized via the restricted functionals \[\mathcal{C}^{\nabla} [1]\Big|_{\widetilde{\Omega}^1(M,T_M)\cap \widetilde{\mathcal{C}}(M)}, \quad \mbox{ respectively } \mathcal{C}^{\nabla}[1, 3]\Big|_{\widetilde{\mathcal{C}}(M)}.\] Moreover, if $(M, g)$ is a Riemannian manifold with Levi-Civita connection $\nabla,$ then based on the discussion at the end of section \ref{csc}, if \[\widetilde{AC}(M)_g:=\{\gamma \in \Omega^{\bullet} (M, T_M) \mid \gamma_1 \in AC(M)_g\},\] then \[\mathcal{C}^{\nabla} [1]\Big|_{\widetilde{\Omega}^1(M,T_M)\cap \widetilde{AC}(M)_g} \quad \mbox{ and } \mathcal{C}^{\nabla}[1, 3]\Big|_{\widetilde{AC}(M)_g}\] realize the K{\"a}hler for $g$ complex structures.
\end{corollary}

The notation $\| \cdot \|$ will denote the norm w.r.t.\ the inner product $\langle \langle \cdot, \cdot \rangle \rangle$ that was introduced earlier. In addition, a top degree form $\alpha \in \Omega^n (M, T_M),$ which can always be expressed in local coordinates $(x_i)_{i=1}^n$ as $\alpha=\sum_{i=1}^n a_i \otimes \frac{\partial}{\partial x_i},$ can be integrated by the rule \[\int_M \alpha:=\sum_{i=1}^n \int_M a_i g\big(\frac{\partial}{\partial x_i}, \frac{\partial}{\partial x_i}\big).\] The same method applies for integrating polyvector bundle top forms in $\Omega^n \big(M, \bigwedge^p {T_M}\big),$ but $g$ should be replaced by the inner product that it induces on $\bigwedge^p {T_M}.$

The formulas from section \ref{int} have the following fact as an immediate consequence.

\begin{proposition}\label{invC}
Let $(M,g)$ be a compact Riemannian manifold of dimension at least $4$ with Levi-Civita connection $\nabla.$ For all $0 \leq k \leq n-2,$ set

\begin{equation*}
\begin{split}
\mathcal{J}(k, \nabla)(A)&:=\|(d^{\nabla})^k \big( d^{\nabla} A \wedge (A \wedge A)-d^{\nabla} A\big)\|^2 \\
&=| (d^{\nabla})^k \big( d^{\nabla} A \wedge (A \wedge A)-d^{\nabla} A\big)|^2_{k+2}.
\end{split}
\end{equation*}
\end{proposition}

Then, \[\mathcal{J}(k,\nabla)=\inf_{A \in AC(M)} \mathcal{J}(k, \nabla)(A)\] is a global numerical obstruction to $M$ being a complex manifold. Moreover, if \[r(\nabla)(A):=\Big| \int_M (d^{\nabla})^{n-2} \big(d^{\nabla} A \wedge (A \wedge A)-d^{\nabla} A\big)\Big|,\] then so is

\[r(\nabla)=\inf_{A \in AC(M)} r(\nabla)(A).\]

For example, if $\dim_{\mathbb{R}} M=6,$ then

\begin{equation*}
\begin{split}
r(\nabla)&=\inf_{A \in AC(M)} \Big| \int_M \big((R^{\nabla})^2 \wedge d^{\nabla} A\big)\wedge (A \wedge A)+8\big((R^{\nabla})^2 \wedge A\big)\wedge (d^{\nabla} A \wedge A)+\\
&12(R^{\nabla} \wedge d^{\nabla} A)\wedge \big((R^{\nabla} \wedge A)\wedge A\big)+8(R^{\nabla} \wedge A) \wedge \big((R^{\nabla} \wedge d^{\nabla} A)\wedge  A\big)\\
&-8(R^{\nabla} \wedge A)\wedge \big((R^{\nabla} \wedge A)\wedge d^{\nabla} A\big)+2d^{\nabla} A\wedge \Big(\big((R^{\nabla})^2 \wedge A\big)\wedge A \Big)+\\
&2d^{\nabla} A\wedge \big((R^{\nabla} \wedge A)\wedge (R^{\nabla} \wedge A)\big)-(R^{\nabla})^2 \wedge d^{\nabla} A \Big|.
\end{split}
\end{equation*}

\begin{proof}
Compactness is tacitly at play in the integrals. The first 2 numerical obstructions are employing Lemma \ref{Q1}, and they can certainly be expanded further using Proposition \ref{PPP1}. The last obstruction is employing the 4th covariant exterior derivative of the intergability form (cf.\ Remark \ref{REM1}). Indeed, if $\mathcal{J}(k,\nabla) >0,$ then $M$ cannot be complex, and the same conclusion can be drawn for $r(\nabla).$
\end{proof}

One might ask what is the advantage of this over working with, for example, the infimum of $\int_M \|N_A\|^2_g \vol_g.$ The answer is that the curvature scalars $\mathcal{J}(k,\nabla),$ and $r(\nabla)$ involve mainly the curvature of $(M,g),$ a quantity that is overall better understood than the Nijenhuis tensor. Consider, for instance, the round $6$-sphere $S^6.$ 

\section*{Acknowledgement}
The author would like to thank the referees for the helpful suggestions. Although this note was written in Paris in an autonomous way, the current version was produced while being supported by the DMS Fellowship of Masaryk University, and the grant project GACR 22-15012J.

\noindent 
G.\ Clemente

\noindent
Department of Mathematics and Statistics, Masaryk University, Building 8, Kotlarska 2, 61137 Brno, Czech Republic

\noindent
e-mail: clemente6171@gmail.com
\end{document}